\documentclass[11pt]{article}
\usepackage{amsmath,amssymb}
\usepackage{graphicx}
\usepackage{hyperref}
\usepackage{authblk}

\newtheorem{thm}{Theorem}
\newtheorem{prp}{Proposition}

\newtheorem{rem}{Remark}

\usepackage{breqn}

\title{Bicentric configurations of pentagonal linkages}
\author[]{A.Diakvnishvili}
\affil[]{Faculty of Business,Technology and Education, Ilia State University, Tbilisi, Georgia}
\date{\today}

\begin{document}

	\maketitle

\begin{abstract}
Let $P$ be a planar $n$-gon with the sidelengths $a_1, \ldots, a_n$ and let us denote by $L=L(P)$ the corresponding planar polygonal linkage. We are concerned with the problem of finding conditions on the sidelengths $a_i$ which guarantee the existence of a bicentric configuration of $L$. Specifically, we present some related results on tangential and chordal polygons and characterize pentagonal linkages having bicentric configurations.
\end{abstract}

		\vspace*{0.25cm}\noindent{\small {\bf Keywords and phrases}:
			{tangential polygon, chordal polygon, bicentric polygon, polygonal linkage, oriented area, inradius, generalized Heron polynomial}
			
			\vspace*{0.25cm}\noindent{\small {\bf MSC 2010:} {52C35,  32S40.}}

			\section*{Introduction}
			
		The aim of this paper is to find explicit conditions on the sidelengths of a planar pentagonal linkage $L$ which guarantee the existence of a bicentric configuration of $L$. To this end we present several general results on tangential and chordal polygons which enable us to give an effective characterization of pentagonal linkages having bicentric configurations. 
		
		Bicentric polygons have been studied in many papers (see, e.g., \cite{jos1}, \cite{josdal}, \cite{min}), in particular, in relation with certain extremal problems (see, e.g., \cite{kea}) and geometric porisms \cite{ber}. As is known from the classical geometry, a polygon with more than three sides in general is not bicentric. This only happens under special conditions on the shape of polygon. For a bicentric polygon $P$, there are important relations between the circumradius $R(P)$, inradius $r(P)$,  and the distance $d(P)$ between the centers of incicrle and circumcircle. The triple $(R(P), r(P), d(P))$ is called the {\it Euler triple} of $P$ and the arising relations between $R(P)$, $r(P)$ and $d(P)$ may be considered as necessary conditions of bicentricity. For triangles and bicentric quadrilaterals, this issue is classical (see, e.g., \cite{jos1}, \cite{josdal}).
		
		For a triangle $ T $ such a relation was given by L.Euler. Namely, we have
		$$
		R^{2}-d^{2}=2Rr,
		$$
		where $ (R, r, d) $  is the Euler triple of $ T $ introduced above.
		
		The first researcher who was concerned with general bicentric polygons was the German mathematician Nicolaus Fuss. He found analogous relations (necessary and sufficient conditions on the triple  $ (R, r, d) $) for bicentric quadrilaterals, pentagons, hexagons, heptagons and octagons. For this reason, those relations are known as Fuss's relations, also in the cases where $ n > 8 $. Nowaday, such relations are known for all $n\leq 20$ \cite{rad2}.
		
		A very remarkable theorem concerning bicentric polygons was given by the French mathematician Victor Poncelet. This theorem can be stated as follows \cite{ber}. If there exists a bicentric $n$-gon with the circumcircle $ C_1 $ and incircle $ C_2 $, then there are infinitely many bicentric $n$-gons whose circumcircle is $ C_1 $ and incircle is $ C_2 $.
		Moreover, such a polygon may any point of $C_1$ ar a vertex (see, e.g., \cite{rad1}). 
		This famous theorem dates of nineteenth century. Since then many mathematicians worked on a number of problems connected with this theorem and solved many of them. 
		
		However there only exist a few papers concerned with bicentric configurations of polygonal linkages (see \cite{bibikhim}, \cite{jos2} and references therein). One of the problems studied in the mentioned papers was concerned with the existence of tangential and bicentric $n$-gon with the prescribed lengths of the sides, which can be naturally interpreted in terms of configurations of polygonal linkages. To the best of our knowledge the problem of existence of bicentric configurations was only solved for $n=4$ (see, e.g., \cite{jos1}, \cite{josdal}, \cite{min}). 
		
		In this paper, we give an effective solution of this problem for $n=5$ and outline possible generalizations of our results. We use an algebraic approach based on the explicit formulas for the inradius, circumradius and oriented area of bicentric polygon, which can be found in \cite{rob1}, \cite{rad1}, \cite{busi}. To make our considerations more clear we present those formulas separately for the tangential and chordal polygons. 
		
		To place our discussion in a proper context we proceed by describing the problem in the general case. Let $P$ be a planar $n$-gon with the positive sidelengths $a_1, \ldots, a_n$. Let us denote by $L=L(P)$ the corresponding planar polygonal linkage. We aim at finding conditions on the sidelengths $a_i$ which guarantee the existence of a bicentric configuration of $L$. 
		
		As is well known, each planar polygonal linkage has a unique convex chordal (aka cyclic) configuration for which the oriented area attains the absolute maximum \cite{khim}, \cite{khpasi}. If this chordal configuration is a tangential polygon then it is bicentric and so $L$ has a bicentric configuration. The uniqueness of convex chordal configuration implies that a convex bicentric configuration may be only one. As was already mentioned, in the case of quadrilateral linkage, criteria of existence of bicentric configuration have been discussed in many papers (see, e.g., \cite{jos1}, \cite{pra}) and nowadays are well known. 
		
		It should be noted that, for $n\geq 5$, the situation becomes essentially complicated since there may also exist non-convex tangential and bicentric configurations. An evident example is given by the regular pentagonal linkage which has two bicentric configurations: convex regular pentagon and regular $5$-pointed star. So for $n\geq 5$ it is natural to consider all possible bicentric configurations including the non-convex ones. We are doing so and this is an essential novelty of the present paper as compared with the aforementioned papers on bicentric configurations of quadrilateral linkages. 
		
		Our main results are as follows. First, we give an effectively verifiable necessary condition for the existence of a bicentric configuration of pentagonal linkage (Theorem \ref{pentabicent1}). Second, we give a characterization of the pentagonal linkages having a convex bicentric configuration (Theorem \ref{pentabicent2}). We also present several illustrative examples and outline a few feasible generalizations of our approach and results.
		
		Our approach relies on the algebraic equations for the inradius of tangential polygon and for the oriented area of chordal and tangential polygons, which we call the modified Sylvester equation (MSE) and modified Robbins equation (MRE) respectively. The original Robbins equation (\ref{RE}) was given in \cite{rob1}, while the modified Sylvester equation (\ref{pentainradius}) is a particular case of the classical results of J.Sylvester  on the inradius of tangential polygon (cf., e.g., \cite{ber}).		   
		
		\section*{\textbf{Background on tangential and chordal polygons}}
		
		We begin with some basic definitions and notations. In this paper, we only deal with the so-called planar polygons lying in a given Euclidean plane endowed with a fixed Cartesian coordinate system.\\
		
		\textbf{Definition 1}. A planar polygon which has a circumscribed circle, i.e. all of its vertices belong to a certain circle, is called a cyclic (or chordal) polygon.
		
		The term "cyclic" is more popular but we stick to the term "chordal" since "cyclic" is overloaded. \\
		
		\textbf{Definition 2}. A planar polygon which has an inscribed circle, i.e. a circle that is tangent to each side of the polygon at its inner point, is called a tangential (or circumscribing) polygon. 
		
		The inscribed circle is usually referred to as the incircle. \\
		
		\textbf{Definition 3}. A planar polygon which simultaneously has a circumscribed circle and an inscribed circle is called a bicentric polygon. \\
		
		In other words, a polygon $P$ is bicentric if it is simultaneously chordal and tangential. For a bicentric polygon, we always denote by $ R=R(P) $ (circumradius) and $ r=r(P) $ (inradius) the radii of circumscribed and inscribed circles respectively. The distance between the centers of these two circles will be called the {\it deflection} $d(P)$ of bicentric polygon $P$. If a bicentric polygon $ P $ is fixed and no misunderstanding can arise, we write simply $ (R, r, d) $ and call it the {\it Euler triple} of $ P $. Obviously, each triangle is bicentric, and each regular polygon is bicentric with the deflection equal to zero. We emphasize that convexity is not assumed in the above definitions so a star-shaped polygon can also be bicentric. 
		
		We will deal with the planar configurations of a given polygonal linkage. Under a polygonal linkage we understand, as usual, the collection of polygons with the fixed lengths of the sides and changeable angles (see, e.g., \cite{ber}). In the next two sections we present several general results on tangential and cyclic polygons which are crucial for our further considerations.  
		
		\section*{Tangent lengths and inradius of tangential polygon}
		
		Many problems concerned with tangential polygons can be conveniently expressed in terms of the so-called circulant matrices \cite{josdal}. In particular, it refers to the existence of a tangential configuration of $n$-gonal linkage with the given sidelengths $a_i$.  Recall that for a tangential polygon $P$ with the incircle $\Gamma$, the {\it tangent length} $t_j$ is defined as the distance between the $j$-th vertex $V_J$ of $P$ and the point $T_j$ at which $\Gamma$ is tangent to the side $[V_jV_{j+1}]$.
		
		Now, assuming that a tangential configuration exists one notices that its tangent lengths should satisfy the following system of linear equations
		\begin{equation} \label{circulant}
			t_j + t_{j+1} = a_j, j=1, \ldots ,n,
		\end{equation}
		where $a_j$ are the lengths of the sides of $P$ (i.e., the sidelengths of the linkage $L(P).$. 
		
		It is known that a tangential configurations of $L(P)$ exists if and only if this system has a solution consisting of positive numbers $t_j$ which we call a positive solution of system (\ref{circulant}) (see, e.g., \cite{busi}). Clearly, the matrix of this system is a circulant $(n\times n)$-matrix of the type
		$$Circ(1,1,0, \ldots, 0).$$
		
		It is well known and easy to verify that the determinant of such matrix is non-zero if and only if $n$ is odd, which implies that the problem of existence of tangential configuration has essentially different solutions for odd $n$ and even $n$. For odd $n$, an $n$-gon linkage can have not more than one tangential configuration. For even $n$, if a tangential configuration exists then there is an infinite set of tangential configurations homeomorphic to a segment \cite{rad2}. 
		
		To make the exposition self-contained we discuss now in some detail the case of tangential pentagon and provide a proof of the following result which just a particular case of an analogous general result given, e.g., in \cite{busi}. We present here an outline of its proof since it yields the formulas for the tangent lengths needed in the sequel.
		
		\begin{prp} \label{circulantmatrix}  	
			If for a given pentagonal linkage $L$ with the sidelengths $a_1, a_2, a_3, a_4, a_5 >0$ the linear system 
			\begin{equation} \label{circulant5}
				\begin{aligned}
					\begin{cases}
						x_1 + x_2 &= a_1 \\
						x_2 + x_3 &= a_2 \\
						x_3 + x_4 &= a_3 \\
						x_4 + x_5 &= a_4 \\
						x_5 + x_1 &= a_5\\
					\end{cases}
				\end{aligned}
			\end{equation}
			has a positive solution, then this solution is unique and the given pentagonal linkage has a tangential configuration. Moreover, a positive solution exists if an only if the sidelengths $a_j$ satisfy a system of explicitly given inequalities
		\end{prp}

		\textbf{Proof}. Notice first that this system always has a non-zero solution. To show this, we calculate the determinant of the system's matrix $A$ and solve it for $x_i,  i=1,2,...,5,$
		where
		$$
		\begin{aligned}
			A = \begin{bmatrix}
				1 & 1 & 0 & 0 & 0 \\
				0 & 1 & 1 & 0 & 0 \\
				0 & 0 & 1 & 1 & 0 \\
				0 & 0 & 0 & 1 & 1 \\
				1 & 0 & 0 & 0 & 1 \\
			\end{bmatrix}
		\end{aligned}
		$$
		
		It is easy to verify that $det(A)=2$ , so this system has a non-zero solution, which is unique. It remains to find conditions which guarantee that this solution is positive. To this end we explicitly solve the system (\ref{circulant5}) for $x_i$, where
		
		\[
		X = \begin{bmatrix}
			x_1 \\
			x_2 \\
			x_3 \\
			x_4 \\
			x_5
		\end{bmatrix}, \quad
		A = \begin{bmatrix}
			1 & 1 & 0 & 0 & 0 \\
			0 & 1 & 1 & 0 & 0 \\
			0 & 0 & 1 & 1 & 0 \\
			0 & 0 & 0 & 1 & 1 \\
			1 & 0 & 0 & 0 & 1
		\end{bmatrix}, \quad
		B = \begin{bmatrix}
			a_1 \\
			a_2 \\
			a_3 \\
			a_4 \\
			a_5
		\end{bmatrix}
		\]
		
		We have $X = A^{-1}B$, where
		
		\[
		A^{-1} = \frac{1}{2} \begin{bmatrix}
			1 & -1 & 1 & -1 & 1 \\
			1 & 1 & -1 & 1 & -1 \\
			-1 & 1 & 1 & -1 & 1 \\
			1 & -1 & 1 & 1 & -1 \\
			-1 & 1 & -1 & 1 & 1
		\end{bmatrix}
		\]
		
		Finally, we get $x_i$ in the form
		
		\begin{align*}
			x_1 &= \frac{1}{2}(a_1 - a_2 + a_3 - a_4 + a_5), \\
			x_2 &= \frac{1}{2}(a_1 + a_2 - a_3 + a_4 - a_5), \\
			x_3 &= \frac{1}{2}(-a_1 + a_2 + a_3 - a_4 + a_5), \\
			x_4 &= \frac{1}{2}(a_1 - a_2 + a_3 + a_4 - a_5), \\
			x_5 &= \frac{1}{2}(-a_1 + a_2 - a_3 + a_4 + a_5).
		\end{align*}
		
		It is now clear that a positive solution exists if and only if all right-hand-sides of the above formulas are strictly positive. The proof of the proposition is complete.\\ 
		
		Notice that we obtained an explicit criterion for the existence of tangential configuration for pentagonal linkage. In the sequel we show that this result yields an effective necessary condition for the existence of a bicentric configuration of pentagon linkage. To this end we will use an important relation between the inradius and tangent lengths. The mentioned relation follows from some classical results of J.Sylvester  (see, e.g., \cite{ber}) but we need its certain modification which will be called the modified Sylvester equation.
		
		\begin{rem}
			It is known that for odd $n$, computation of the inradius can be reduced  to finding the positive real roots of a certain algebraic equation. At the same time, it is also known that, for even $n$, the possible values of inradii of an $n$-gonal linkage $L$ which has tangential configurations form a certain segment. In this case, a natural problem is to find the endpoints of this segment, i.e., the minimal and maximal possible values of the inradii of tangential configurations of $L$. For $n=4$ this easily follows from the known results but the case of even $n\geq 6$ requires more comprehensive analysis beyond the scopes of this paper. This is the main reason why we restrict our considerations to pentagonal linkages.
		\end{rem}
		
		As we have seen, the tangent lengths of a tangential pentagon are expressible as linear functions of its sidelengths. At the same time, as was already said the inradius of tangential polygon can be found from the tangent lengths by solving an algebraic equation. This follows from some classical results going back to J.Sylvester \cite{ber}. We present a general result and then specify if for $n=5$.  
		
		Let us again denote by $t_j$ the tangent lengths of an $n$-gon $P$ and by $s_k$ the $k$-th elementary symmetric functions of $t_j$. Next, introduce a polynomial $I_n(x)$ in one variable $x$ using the following recipe:\\
		If $n=2k+1$ then $I_n(x) = s_1x^k - s_3x^{k-1} + \ldots + (-1)^ks_n.$ \\
		If $n=2k+2$ then $I_n(x) = s_1x^k - s_3x^{k-1} + \ldots + (-1)^ks_{n-1}$.
		
		It is known that in both these cases the polynomial $I_n$ has exactly $k$ positive real roots which are equal to the squares of the inradii of tangential polygons, not necessarily convex, with the tangent lengths $t_j$ \cite{busi}.  
		
		For $n=5$ we get the following conclusion which is crucial for our purposes. The equation given will be called the Sylvester equation for pentagon (see, e.g., \cite{busi}).  A self-contained proof of this result based on the use of complex numbers can be found in \cite{busi}.
		
		\begin{prp} \label{PentainRadius}
			The inradius of a tangential pentagon with tangent lengths $t_i$ can be found as a positive root of the following quartic equation	
			\begin{equation} \label{pentainradius}
				s_1 x^4 - s_3 x^2 + s_5 = 0, 
			\end{equation}
			where $\sigma_j$ are the elementary symmetric polynomials in $t_i$ 
		\end{prp}

		Since the equation (\ref{pentainradius}) is biquadratic we conclude that the two possible values of the inradii can be explicitly computed through the tangent lengths $t_j$ and also through
		the sidelengths $a_j$.

		\section*{Oriented area of chordal polygon}
		
		As is known from the classical geometry, one can compute the area of triangle and cyclic quadrilateral from their sides by the well-known Heron's and Brahmagupta's formulas. Generalizations of these classical results were given by an American mathematician David Robbins \cite{rob1}, \cite{rob2}. More precisely, D.Robbins gave a method of constructing certain polynomials in one variable, called {\it generalized Heron polynomials} (GHP), real roots of which provided information on the oriented area of cyclic polynomials \cite{rob1}. 
		
		In particular, D.Robbins found a polynomial which enables one to calculate the area of a cyclic pentagon through its sides only \cite{rob1}. For this reason the generalized Heron polynomial in the case of pentagon will be referred as {\it Robbins' polynomial.} For clarity and further use we give now a precise formulation of the original result of D.Robbins.  
		
		\textbf{Robbins' polynomial}. Given positive real numbers $a_1, . . . . ., a_5$, we seek for the extremal values of the oriented area of chordal pentagons, convex or not, whose side lengths are $a_1, a_2, a_3, a_4, a_5$. D.Robbins constructed the following polynomial, whose maximal root is equal to the sixteen times square of the convex cyclic pentagon's area. Namely, if we put $u=16\,{A}^{2},$ where $A$ is the area of such a cyclic pentagon, then the unknown $u$ can be found as a real root of the following polynomial usually denoted by $H_5$:
		
		\begin{equation} \label{RE}
			H_5(u) = -27\,{u}^{2}{t_{{5}}}^{2}-18\,ut_{{3}}t_{{4}}t_{{5}}+u{t_{{4}}}^{3}-16
			\,{t_{{3}}}^{3}t_{{5}}+{t_{{3}}}^{2}{t_{{4}}}^{2},
		\end{equation}
		where 
		$$
		t_{{2}}={e_{{1}}}^{2}+u-4\,e_{{2}},
		$$
		$$
		t_{{3}}=e_{{1}}t_{{2}}+8\,e_{{3}}={e_{{1}}}^{3}+ue_{{1}}-4\,e_{{1}}e_{{2}}+8\,e_{{3}},
		$$
		$$
		t_{{4}}={e_{{1}}}^{4}+2\,u{e_{{1}}}^{2}-8\,{e_{{1}}}^{2}e_{{2}}+{u}^{2
		}-8\,ue_{{2}}+16\,{e_{{2}}}^{2}-64\,e_{{4}},
		$$
		$$
		t_{{5}}=128\,e_{{5}},
		$$
		and $e_j, j = 1, 2, .. , 5$ are the elementary symmetric functions in the squares of the sides $a_j$ (see \cite{rob1}).
		
		It is easy to see that if one the sides is equal to zero we obtain an equation equivalent to Brahmagupta formula, which in turn is a generalization of Heron's formula. D.Robbins also outlined a method of constructing generalized Heron polynomials for arbitrary $n$ and formulated several conjectures about the degrees and coefficients of such GHPs. These conjectures have been proved in papers of R.Connelly, M.Fedorchuk and I.Pak (see, e.g., \cite{pak}). We present an extract from their results in the form relevant for our purposes. 
		
		For any $n$-tuple positive real numbers $a = (a_1, . . . . ., a_n),$ there exists a monic polynomial $H_n(a)$ in one variable of degree $d_n$ such that the oriented area of any planar cyclic $n$-gon with the sidelengths $a_i$ is equal to one of the roots of $H_n(a)$. The coefficients of $H_n(a)$ are symmetric functions of the squares $a_i^2$ of the given sidelengths and its degree $d_n$ can be calculated as follows. Put
		$$\triangle_k = \frac{(2k+1)\times 2k \times ... \times (k+1)}{2(k!)}.$$
		Then if $n=2k+1$ we have $d_n = \triangle_k$, and if $n=2k+2$ we have $d_n = 2\triangle_k.$ 
		In particular, for $n=5$ we have $d_5 = 7,$ which is equal to the degree of the Robbins' polynomial \ref{RE}. It turns our that polynomial $H_5$ can be used in the study of existence problem for bicentric configurations of pentagonal linkages.

		\section*{Existence problem for bicentric configurations}
		
		In this section we show how the results presented above can be used to solve the existence problem for bicentric configurations. Suppose we have a bicentric configuration $X$ of $n$-gonal linkage $L$. Then we have $A = rp$ where $p$ is its semiperimeter. To obtain an effective necessary condition for the existence of bicentric polygon we rewrite the (\ref{pentainradius}) in terms of the area $A$ by substituting $r = A/p$ and clearing the denominators. This gives us the following algebraic equation in $z$ which should be satisfied by the area $A$ of a bicentric pentagon linkage:
		
		\begin{equation} \label{MMSE}
			s_1 z^4 - p^2 s_3 z^2 + p^4 s_5 = 0. 
		\end{equation} 
		
		Notice that the left-hand-side is in fact a polynomial in $z^2$. The modified Sylvester polynomial $\tilde I_5(u)$ is now obtained by passing from $z$ to the variable $u = 16z^2$ which reduces to just rescaling of coefficients. Now notice that, by the construction of polynomial $\tilde I_5(u)$, for a bicentric configuration the value of $u= 16A^2$ should be a common root of $\tilde I_5$ and $H_5$. Remembering that the existence of a common root implies the vanishing of the resultant we arrive at our first main result. 
		
		\begin{thm} \label{pentabicent1}
			For a pentagonal linkage $L$, a necessary condition for the existence of a bicentric configuration is that
			$$Res_u(\tilde I_5, H_5) = 0,$$
			where $Res_u$ denoted the resultant of the two polynomials in variable $u.$
		\end{thm}
		
		This necessary condition is effective since the coefficients of $\tilde I_5$ and $H_5$ are algebraically computable the the sidelengths of $L$. Now we can also give a criterion of existence of a convex bicentric configuration.  
		
		\begin{thm} \label{pentabicent2}
			If we have a tangential convex configuration of pentagonal linkage with the given sidelengths and the sixteen times square of its area is the maximal root of the Robbins' polynomial for the same sides, then this configuration is bicentric.
		\end{thm}
		
		\textbf{Proof.} From the given sidelengths we compute the tangent lengths and coefficients of the modified Sylvester equation. This enables us to explicitly compute the two possible values of the oriented area $A$ of tangential configuration and insert them in the Robbin's polynomial. Since it is known that the convex cyclic configuration is determined by its area \cite{bibikhim} it follows that if the maximal value of $16A^2$ is a root of the Robbin's polynomial then the corresponding tangential configuration is at the same cyclic. The proof is complete. \\
		
		In the case of non-convex tangential configuration the situation is more complicated and requires additional considerations described below. 
		
		\begin{prp} \label{Pentabicent3}
			The value of the inradius of a non-convex tangential configuration can be explicitly computed form the so-called inverse tangent equation. 
		\end{prp}
		
		\textbf{Proof.} From the given sidelengths we compute the solution $(t_i)$ of the linear system for the tangent lengths This enables us to compute the inradius $r_1$ of a non-convex tangential configuration from the so-called inverse tangent equation which for non-convex configuration has the form
		
		$$
		\sum_{i=1}^{5} \arctan\left(\frac{t_i}{r_1}\right) = \pi.
		$$
		
		This in turn enables us to compute the value $16A^2$ for the non-convex tangential configuration and insert it in the Robbin's polynomial to check if it is cyclic. In this way one can also give a sufficient condition for the existence of non-convex bicentric configuration.
		
		The case of non-convex tangential configuration has some interesting aspects concerned with the so-called Whitney index of tangential polygon (see, e.g., \cite{busi}) but we do not discuss them here for the reason of space. Instead we present for completeness the analogous results for quadrilateral linkages. 
		
		\section*{Bicentric configurations of quadrilateral linkage}
		
		As was mentioned, the case of quadrilateral linkage was well studied (see, e.g., \cite{jos1} or \cite{min}). For this reason we present the corresponding criterion, usually called Pitot's criterion, without proof.
		
		\begin{prp} (Pitot's criterion) \label{quadrocriterion}
			A nondegenerate quadrilateral linkage $L(a_1, \ldots, a_n)$ with all $a_i>0$ has a bicentric configuration 
			if and only if
			
			$$a_1 + a_3 = a_2 + a_4.$$
		\end{prp}
		
		It follows that a quadrilateral linkage satisfying the Pitot's criterion has a whole continuum of tangential configurations. By the way of analogy with the extremal problems for open polygonal linkages considered in \cite{bibikhim} it is natural to wonder what are the extremal values of inradii of tangential configurations of a quadrilateral linkage satisfying the Pitot's criterion. 
		
		The maximal value $r^*$ is easily available due to the maximality of the cyclic configuration \cite{khim}, Brahmagupta formula, and evident relation between the area and perimeter of a tangential polygon, which yields the following result. 
		
		\begin{prp}
			$r^* = \frac{\sqrt{(p-a)(p-b)(p-c)(p-d)}}{p}.$
		\end{prp}
		
		To find the minimal value of inradius we need to consider two cases. By Pitot's criterion such a linkage has aligned configurations. If it is a kite then it is easy to show that the minimal value of inradiii is zero. If it is not a kite then it has two triangular configurations with two aligned sides and an elementary geometric argument shows that the minimal value of inradii is the inradius of one of the arising triangles. Thus for quadrilateral linkages all the problems formulated above have simple 
		and effective solutions.

		\section*{\textbf{Examples}} 
		
		We now illustrate the main results by calculating the areas of tangential and cyclic pentagons in several examples.\\
		
		\textbf{Example 1. Area of convex tangential regular pentagon}\\
		
		If we have a regular pentagon with the sides $a_1=a_2=...=a_5=1$, then its area $A$ is given by the well-known following formula \\
		$$
		2A=Pr,
		$$
		where $P$ is the perimeter and $r$ is the inradius. So we get that \\
		$P = 5$, $r=0.6881909602$.\\
		Hence the area is equal to \\
		$A = 1.720477400$, which coincides with the maximal root of Robbins' polynomial.\\
		
		\textbf{Example 2. Area of convex tangential non-regular pentagon}\\
		
		If we have a pentagon with different sides $a_1 = 29, a_2 = 30, a_3 = 31, a_4 = 32, a_5 = 33$, then its area (denoted as $A_1$) is calculated by the same formula \\
		$$
		2A_1=P_1r_1,
		$$
		where 
		$P_1$ and  $r_1$ are, respectively, the perimeter and inradius. From the given sides we get that \\
		$P_1= 155$, where \\
		$r_1$ is the solution of the following equation
		$$
		\sum_{i=1}^{5} \arctan\left(\frac{x_i}{r_1}\right) = \pi.
		$$
		Here $x_i$ is calculated from the linear system $(0.2)$. In this example we have equation\\
		$$
		\arctan\left(\frac{13.5}{r_1}\right) + \arctan\left(\frac{14.5}{r_1}\right) + \arctan\left(\frac{15.5}{r_1}\right) + \arctan\left(\frac{16.5}{r_1}\right) + \arctan\left(\frac{17.5}{r_1}\right) = \pi,
		$$
		from which we get $r_1=21.27248379$, hence $A_1=1648.617494$. As is shown below in this case
		this is again the maximal root of the Robbin's polynomial. \\
		
		\textbf{Example 3. Area of a convex cyclic regular pentagon} \\
		
		If we have a cyclic pentagon with sides $a_1=a_2=...=a_5=1$ and want to calculate its area $K$ , first write an elementary symmetric functions for Robbins' polynomial \\
		$$
		e_1 = 5, e_2 = 10, e_3 = 10, e_4 = 5, e_5 = 1.
		$$
		Now, if we  insert $e_i  , i=1,2,...,5$ in Robbin's polynomial, we  get the  following  equation \\
		$$
		{u}^{7}-65\,{u}^{6}+965\,{u}^{5}-6645\,{u}^{4}+25155\,{u}^{3}-54243\,{
			u}^{2}+62775\,u-30375=0
		$$
		with the roots: $2.639320225, 3., 3., 3., 3., 3., 47.36067977.$ \\
		So the maximal root is $47.36067977$ and we get
		$$
		16 K^2 = 47.36067977
		$$
		$K =1.720477400$. \\
		
		\textbf{Example 4. Area of convex cyclic irregular pentagon} \\
		
		Now we use the same argument for a irregular pentagon with pairwise different sides $a_1=29, a_2=30, a_3=31, a_4=32,a_5=33$. An area denotes by $K_1$, we write the elementary symmetric functions for Robbins' polynomial \\
		\begin{equation}
			\begin{aligned}
				e_1 &= 4815 \\
				e_2 &= 9254463 \\
				e_3 &= 8875070485 \\
				e_4 &= 4246737436836 \\
				e_5 &= 811128627302400
			\end{aligned}
		\end{equation}
		If we  insert $e_i, i=1,2,...,5$ in Robbins' polynomial,  we  get the  following equation
		$$
		\begin{aligned}
			{u}^{7}-59817537\,{u}^{6}+814568856314373\,{u}^{5}- \\
			5129732167330152025589\,{u}^{4}+17708992633706617259476903875\,{u}^{3}\\
			-34729928934462676267203902962651875\,{u}^{2}+ \\
			36459248759033130575200748650105233984375\,u- \\
			15963113698969945651994827119257850429052734375=0
		\end{aligned}
		$$
		with the roots: $2.337566549\cdot 10^6, 2.350203962\cdot 10^6, 2.499491552\cdot 10^6, 2.713188068\cdot 10^6, 2.979318223\cdot 10^6, 3.295261735\cdot 10^6, 4.364250691\cdot 10^7$. \\
		
		Thus the maximal root is  $4.364250691\cdot 10^7=43642506.91$ and find that
		$$
		16\,{K_{{1}}}^{2}= 43642506.91.
		$$
		Hence $K_1=1651.561892$.\\
		
		\section*{Concluding remarks}
		
		There are a number of open problems and research prospects suggested by the results of this paper. Below me mention some of them in a hope that they may attract attention of other researchers and lead ot further progress in this field.
		
		One of the most natural problems is to extend our approach to linkages with the number of sides $n\geq 6$. For odd values of $n$ one may use the same approach as above. In particular, for $n = 7, 9$ one can get a criterion analogous to Theorem \ref{bicent2} because, for these values of $n$, the degree of the modified Sylvester equation does not 
		exceed $4$. Thus the roots of the latter equation are expressible through the sidelengths and can be inserted in the generalized Heron equations.
		
		However, for even values of $n$ one needs additional considerations cause by the existence of continua of tangential configurations. In particular, a topical problem is to find the minimal and maximal possible values of the inradii in this case. 
		
		It is also natural and feasible to extend our approach and results to the open polygonal linkages also known as robot $n$-arms. For $n\geq 3$, the bicentric configurations may form continual subsets in the multidimensional configuration space and there arises the problem of characterization of their possible topological and geometric structure.
		First few steps in this direction were made in \cite{bibikhim} but only for $n\leq 4$.
		
		Finally, the same problems are meaningful and feasible for spherical linkages and linakes in the hyperbolic plane. To the best of our knowledge these topics remain practically unexplored and our methods may lead to progress in these fields.

		\bigskip
		
		Author's e-mail and address:\\
		
		Ana Diakvnishvili\\
		Ilia State University\\
		Faculty of Business, Technology and Education\\
		3/5, K. Cholokashvili Ave., Tbilisi 0162. Georgia\\
		E-mail: ana.diakvnishvili.1@iliauni.edu.ge\\

	\end{document}